\documentclass[11pt]{article}
\usepackage{amssymb,amsfonts,amsmath,latexsym,epsf,tikz,url}

\newtheorem{theorem}{Theorem}[section]

\newtheorem{corollary}[theorem]{Corollary}
\newtheorem{lemma}[theorem]{Lemma}

\newcommand{\proof}{\noindent{\bf Proof.\ }}
\newcommand{\qed}{\hfill $\square$\medskip}

\begin{document}

\title{The distinguishing number and the distinguishing index of  co-normal product of two graphs}

\author{
	Saeid Alikhani
\and 
 Samaneh Soltani$^{}$\footnote{Corresponding author}
}

\date{\today}

\maketitle

\begin{center}
Department of Mathematics, Yazd University, 89195-741, Yazd, Iran\\
{\tt alikhani@yazd.ac.ir, s.soltani1979@gmail.com}
\end{center}


\begin{abstract}
The distinguishing number (index) $D(G)$ ($D'(G)$)  of a graph $G$ is the least integer $d$
such that $G$ has an vertex labeling (edge labeling)  with $d$ labels  that is preserved only by a trivial automorphism. The co-normal product $G\star H$ of two graphs $G$ and $H$ is the graph with vertex set $V (G)\times V (H)$ and edge set $\{\{(x_1, x_2), (y_1, y_2)\} | x_1y_1 \in E(G) ~{\rm or}~x_2y_2 \in E(H)\}$.
In this paper we study the distinguishing number and the distinguishing index of the co-normal product of two graphs.  We prove that for every $k \geq 3$, the $k$-th co-normal power of a connected 
graph $G$ with no false twin vertex and no dominating vertex, has the distinguishing number and the distinguishing  index equal two.

\end{abstract}


\medskip
\noindent{\bf AMS Subj.\ Class.:} 05C15, 05C60. 

\section{Introduction and definitions}
Let $G=(V,E)$ be a simple graph of order $n\geq 2$. We use the the following notations: The set of vertices adjacent in $G$ to a vertex of a vertex subset $W\subseteq  V$ is the \textit{open neighborhood}  $N(W )$ of $W$. Also $N(W )\cup W$ is called a \textit{closed neighborhood}   of $W$ and denoted by $N[W]$. 
   A  \textit{subgraph} of a graph $G$ is a graph $H$ such that 
$V(H) \subseteq V(G)$ and  $E(H) \subseteq E(G)$. If $V(H) = V(G)$, we call $H$ a \textit{spanning subgraph} of $G$. Any spanning 
subgraph of $G$ can be obtained by deleting some of the edges from $G$. Two distinct vertices $u$ and $v$ are called \textit{true twins} if $N[v] = N[u]$ and \textit{false twins}  if $N(v) = N(u)$. Two vertices are called \textit{twins} if they are true or false twins. The
number $|N(v)|$ is called the \textit{degree} of $v$ in $G$, denoted as ${\rm deg}_G(v)$ or ${\rm deg}(v)$. A vertex having degree $|V(G)| - 1$ is called a \textit{dominating vertex}  of $G$. 
 Also, ${\rm Aut}(G)$ denotes the automorphism group of $G$, and graphs with  $|{\rm Aut}(G)|=1$ is called \textit{rigid} graphs.

A labeling of $G$, $\phi : V \rightarrow \{1, 2, \ldots , r\}$, is said to be \textit{$r$-distinguishing}, 
if no non-trivial  automorphism of $G$ preserves all of the vertex labels.
The point of the labels on the vertices is to destroy the symmetries of the
graph, that is, to make the automorphism group of the labeled graph trivial.
Formally, $\phi$ is $r$-distinguishing if for every non-trivial $\sigma \in {\rm Aut}(G)$, there
exists $x$ in $V$ such that $\phi(x) \neq \phi(\sigma(x))$. The \textit{distinguishing number} of a graph $G$ is defined  by
\begin{equation*}
D(G) = {\rm min}\{r \vert ~ G ~\textsl{\rm{has a labeling that is $r$-distinguishing}}\}.
\end{equation*} 

This number has defined in \cite{Albert}. Similar to this definition,  the \textit{distinguishing index}  $D'(G)$ of $G$ has defined in \cite{Kali1} which is  the least integer $d$
such that $G$ has an edge colouring   with $d$ colours that is preserved only by a trivial
automorphism. If a graph has no nontrivial automorphisms, its distinguishing number is  $1$. In other words, $D(G) = 1$ for the asymmetric graphs.
The other extreme, $D(G) = \vert V(G) \vert$, occurs if and only if $G$ is a complete graph. The distinguishing index of some examples of graphs was exhibited in \cite{Kali1}. For 
instance, $D(P_n) = D'(P_n)=2$ for every $n\geq 3$, and 
$D(C_n) = D'(C_n)=3$ for $n =3,4,5$,  $D(C_n) = D'(C_n)=2$ for $n \geq 6$.  A graph and its complement, always have the same automorphism group while their graph structure usually differs, hence $D(G) = D(\overline{G})$ for every simple graph $G$.

 Product graph of two graphs $G$ and $H$ is a new graph having the vertex set
$V (G) \times V (H)$ and the adjacency of vertices is defined under some rule using the
adjacency and the nonadjacency relations of $G$ and $H$. 
The distinguishing number and the distinguishing index of some graph products has been studied in literature (see \cite{Imrich,Imrich and  Klavzar}). 
The \textit{Cartesian product} of graphs $G$ and $H$ is a graph, denoted by $G\Box H$, whose vertex
set is $V (G) \times V (H)$. Two vertices $(g, h)$ and $(g', h')$ are adjacent if either $g = g'$ and
$hh' \in E(H)$, or $gg' \in E(G)$ and $h = h'$.  In 1962, Ore \cite{ore} introduced
a product graph, with the name Cartesian sum of graphs. Hammack et al.
\cite{Sandi}, named it co-normal product graph.
 The \textit{co-normal product} of $G$ and $H$ is the graph denoted by $G \star H$, and is defined as follows: 
\begin{align*}
&V (G \star H) = \{(g, h) | g \in V (G) ~{\rm and}~ h \in V (H)\},\\
&E(G \star H) =\{\{(x_1, x_2), (y_1, y_2)\} | x_1y_1 \in E(G) ~{\rm or}~x_2y_2 \in E(H)\}.
\end{align*}

 We need knowledge of the structure of the automorphism group of the Cartesian product, which was determined by Imrich \cite{Imrich1969}, and independently by Miller \cite{Miller}.

\begin{theorem}{\rm \cite{Imrich1969,Miller}}\label{autoCartesian} Suppose $\psi$  is an automorphism of a connected graph $G$ with prime
	factor decomposition $G = G_1\Box G_2\Box \ldots \Box G_r$. Then there is a permutation $\pi$ of the set
	$\{1, 2, \ldots ,  r\}$ and there are isomorphisms  $\psi_i : G_{\pi(i)} \rightarrow G_i$, $i = 1,\ldots , r$, such that
	$$\psi(x_1, x_2 , \ldots , x_r) = ( \psi_1(x_{\pi(1)}),  \psi_2(x_{\pi(2)}), \ldots ,  \psi_r(x_{\pi(r)})).$$ 
\end{theorem}

Imrich and Kla\v{v}zar in \cite{Imrich and  Klavzar}, and Gorzkowska et.al. in \cite{A. Gorzkowska R. Kalinowski and M. Pilsniak} showed that the distinguishing number and the distinguishing index of the square and higher powers of a connected graph $G \neq K_2, K_3$ with respect to the Cartesian product is 2.

 The relationship between the automorphism group of co-normal product of two non isomorphic, non rigid connected  graphs with no false twin and no dominating vertex  is the same as that in the case of the Cartesian product.
\begin{theorem}{\rm \cite{fixing}}\label{thm2.12} For any two non isomorphic, non rigid graphs $G$ and $H$, ${\rm Aut}(G\star H) = {\rm Aut}(G) \times {\rm Aut}(H)$ if and only if  both $G$ and  $H$ have no false twins and dominating vertices.
\end{theorem}
\begin{theorem}{\rm \cite{fixing}}\label{thm2.17} For any two rigid isomorphic graphs $G$ and $H$, ${\rm Aut}(G\star H) \cong S_2$.
\end{theorem}
\begin{theorem}{\rm \cite{fixing}}\label{thm2.18}The graph $G\star H$ is rigid if and only if $G \ncong H$ and both $G$ and  $H$ are rigid graphs.
\end{theorem}
\medskip

In the next section, we study the distinguishing number of the co-normal product of two graphs. In section 3, we show that the distinguishing index of the co-normal product of two simple connected non isomorphic, non rigid graphs with no false twin and no dominating vertex  cannot be more than the distinguishing index of their Cartesian product. As a consequence, we prove that all powers of a connected graph $G$ with no false twin and no dominating vertex distinguished by exactly two edge labels  with respect to the co-normal product.

\section{Distinguishing number of co-normal product of two graphs}
We begin this section with a general upper bound for the co-normal product of two simple connected graphs. We need the following theorem.
\begin{theorem}{\rm \cite{fixing}}\label{Theorem 2.4} Let $G$ and $H$ be two graphs and $\lambda : V (G\star H) \rightarrow V (G\star H)$ be a mapping.
\begin{enumerate}
\item[(i)]  If $\lambda = (\alpha, \beta)$ defined as $\lambda(g , h) = (\alpha(g), \beta(h))$, where $\alpha \in {\rm Aut}(G)$ and $\beta \in {\rm Aut}(H)$, then $\lambda$ is an automorphism on $G\star H$.
\item[(ii)] If $G$ is isomorphic to $H$ and $\lambda = (\alpha, \beta)$ defined as $\lambda(g , h) = ( \beta(h),\alpha(g))$
where $\alpha$ is an isomorphism on $G$ to $H$ and $\beta$ is an isomorphism on $H$ to $G$, then
$\lambda$ is an automorphism on $G\star H$.
\end{enumerate}
\end{theorem}
\begin{theorem}\label{generalbound}
If $G$ and $H$ are  two simple connected graphs, then 
\begin{equation*}
{\rm max}\big\{D(G\Box H), D(G), D(H)\big\}\leq D(G\star H)\leq  {\rm min} \big\{D(G)|V(H)|, |V(G)|D(H)\big\}.
\end{equation*}
\end{theorem}
\proof We first show that ${\rm max}\{D(G), D(H)\}\leq D(G\star H)$. By contradiction, we assume that  $D(G\star H)< {\rm max}\{ D(G), D(H)\}$. Without loss of generality we suppose that ${\rm max}\{ D(G), D(H)\}= D(G)$. Let $C$ be a $(D(G\star H))$-distinguishing labeling of $G\star H$. Then the set of vertices $\{(g,h^*)~:~g\in V(G)\}$ where $h^*\in V(H)$ have been labeled with less than $D(G)$ labels. Hence we can define the labeling $C'$ with $C'(g):= C(g,h^*)$ for all $g\in V(G)$. Since $D(G\star H)< D(G)$, so $C'$ is not a distinguishing labeling of $G$, and so there exists a nonidentity automorphism $\alpha$ of $G$ preserving the labeling $C'$. Thus there exists a nonidentity automorphism $\lambda$ of $G\star H$ with $\lambda (g,h):= (\alpha (g), h)$ for   $g\in V(G)$ and $h\in V(H)$, such that $\lambda$ preserves the distinguishing labeling $C$, which is a contradiction. Now we show that $D(G\Box H) \leq D(G\star H)$, and so we prove the left inequality. By Theorems \ref{autoCartesian} and \ref{Theorem 2.4}, we can obtain that ${\rm Aut} (G\Box H)\subseteq {\rm Aut} (G\star H)$, and  since $V(G\Box H)=V(G \star H)$, we have $D(G\Box H) \leq D(G\star H)$. 

Now we show that $ D(G\star H)\leq  {\rm min}\left\{D(G)|V(H)|, |V(G)|D(H)\right\}$. For this purpose, we define two distinguishing labelings of $G\star H$ with $D(G)|V(H)|$ and $|V(G)|D(H)$ labels, respectively. Let $C$ be a $D(G)$-distinguishing labeling of $G$ and $C'$  be a $D(H)$-distinguishing labeling of $H$. We suppose that $V(G)=\{g_1, \ldots , g_n\}$ and $V(H)=\{h_1, \ldots , h_m\}$, and define the two following distinguishing labelings $L_1$ and $L_2$ of $G\star H$ with $D(G)|V(H)|$ and $|V(G)|D(H)$ labels.
\begin{align*}
& L_1(g_j,h_i):= (i-1)D(G)+C(g_j),\\
& L_2(g_j,h_i):= (j-1)D(H)+C'(h_i).
\end{align*}

We only prove that the labeling  $L_1$ is a distinguishing labeling, and by a similar argument, it can be concluded that $L_2$ is a distinguishing labeling of $G\star H$. If $f$ is an automorphism of $G\star H$ preserving the labeling $L_1$, then $f$ maps the set $H_i := \{(g_j,h_i)~:~g_j\in V(G)\}$ to itself, setwise, for all $i = 1, \ldots ,m$. Since the restriction of $f$ to $H_i$ can be considered as an automorphism of $G$ preserving the distinguishing labeling $C$, so for every $1\leq i \leq m$, the restriction of $f$ to $H_i$ is the identity automorphism. Hence $f$ is the identity automorphism of $G\star H$.\qed

The bounds of Theorem \ref{generalbound} are sharp. For the right inequality it is sufficient to consider  the complete graphs as the graphs $G$ and $H$. In fact, if $G= K_n$ and $H= K_m$, then $G\star H = K_{nm}$. For the left  inequality we consider the non isomorphic rigid graphs as the graphs $G$ and $H$. Then by Theorem \ref{thm2.18}, we conclude that $G\star H$ and $G\Box H$ are a rigid graph and hence ${\rm max}\big\{D(G\Box H), D(G), D(H)\big\}= D(G\star H)$.

\medskip
With respect to Theorems \ref{autoCartesian} and \ref{thm2.12}, we have that the automorphism group of a co-normal product of connected  non isomorphic, non rigid graphs with no false twin and no dominating vertex, is the same as automorphism group of the Cartesian product of them, so the following theorem follows immediately:

\begin{theorem}\label{disnumstroncartes}
If $G$ and $H$ are two simple connected, non isomorphic, non rigid  graphs with no false twin and no dominating vertex, then $D(G\star H)=D(G\Box H)$.
\end{theorem}

Since the path graph $P_n$ ($n\geq4$), and the cycle graph $C_m$ ($m\geq5$) are connected,  graphs with no false twin and no dominating vertex, then by Theorem \ref{disnumstroncartes} we have $D(P_n \star P_q)= D(P_n \star C_m)=D(C_m \star C_p)=2$ for any $q,n\geq3$ where $q\neq n$ and $m,p\geq5$ where $m\neq p$. (see \cite{Imrich and  Klavzar} for the distinguishing number of Cartesian product of these graphs).

\medskip
To prove the next result, we need the following lemmas.
\begin{lemma}{\rm \cite{resolving}}\label{Lemma 2.3}
For any two distinct vertices $(v_i,u_j)$ and $(v_r,u_s)$ in $G\star H$, $N((v_i,u_j)) = N((v_r,u_s))$ if and only if 
\begin{enumerate}
\item[(i)] $v_i = v_r$ in $G$ and $N(u_j) = N(u_s)$ in $H$, or 
\item[(ii)] $u_j = u_s$ in $H$ and $N(v_i) = N(v_r)$ in $G$, or 
\item[(iii)] $N(v_i) = N(v_r)$ in $G$ and $N(u_j) = N(u_s)$.
\end{enumerate}
\end{lemma}
\begin{lemma}{\rm \cite{resolving}}\label{Lemma 4.1}
A vertex $(v_i,u_j)$ is a dominating vertex in $G\star H$ if and only if $v_i$ and $u_j$ are
dominating vertices in $G$ and $H$, respectively.
\end{lemma}
\begin{theorem}{\rm \cite{fixing}}\label{thm 2.15} For a rigid graph $G$ and a non rigid graph $H$, $|{\rm Aut}(G\star H)| =|{\rm Aut}(H)|$ if and only if $G$ has no dominating vertex and $H$ has no false twin.
\end{theorem}

Now we are ready to state and prove the main result of this section.
\begin{theorem}
	Let $ G$ be a  connected  graph with no false twin and no dominating vertex,	and $\star G^k$ the $k$-th power of $G$ with respect to the co-normal product. Then $D(\star G^k) = 2$ for $k \geq 3$. In particular, if $G$ is a rigid graph then for $k \geq 2$, $D(\star G^k) = 2$. 
\end{theorem} 
 \proof By Lemmas \ref{Lemma 2.3} and \ref{Lemma 4.1}, we can conclude that $G\star G$ has no false twin and no dominating vertex. We consider  the two following cases:
 
 Case 1)  Let $G$ be a non rigid graph. If $H := G\star G$, then $D(\star G^3) =2$ by Theorem \ref{disnumstroncartes}. Now by induction on $k$, we have the result.
 
 Case 2) Let $G$ be a rigid graph. In this case, $|{\rm Aut}(G\star G)|=2$,  by Theorem \ref{thm2.17}, and so $D(G\star G) =2$. If $H:= G\star G$, then $|{\rm Aut}(G\star H)|=|{\rm Aut}(H)|$, by Theorem \ref{thm 2.15}. Hence $|{\rm Aut}(\star G^3)|=2$. By induction on $k$ and using Theorem \ref{thm 2.15}, we obtain  $D(\star G^k) = 2$  for $k \geq 2$ where $G$ is a rigid graph.\qed

\section{Distinguishing index of co-normal product of two graphs}
In this section we investigate the distinguishing index of co-normal product of graphs.  Pil\'sniak in \cite{nord} showed that the distinguishing index of traceable graphs, graphs with a Hamiltonian path, of order equal or  greater than seven is at most two.
  \begin{theorem}{\rm \cite{nord}}\label{traceble}
 If $G$ is a traceable graph of order $n\geq  7$, then $D'(G) \leq 2$.
 \end{theorem}
 
We say that a graph $G$ is almost spanned by a subgraph $H$ if $G-v$, the graph obtained from $G$ by removal of a vertex $v$ and all edges incident to $v$, is spanned by $H$ for some $v \in  V (G)$. The following two observations will play a crucial role
in this section.
\begin{lemma} {\rm \cite{nord}}\label{nordspanning}
If a graph $G$ is spanned or almost spanned by a subgraph $H$, then $D'(G) \leq D'(H) + 1$.
\end{lemma}
\begin{lemma}\label{distindspann}
Let $G$ be a graph and $H$ be a spanning subgraph of $G$. If ${\rm Aut}(G)$ is a subgroup of ${\rm Aut}(H)$, then $D'(G)\leq  D'(H)$.
\end{lemma}
\proof  Let to  call the edges of $G$ which are the edges of $H$, $H$-edges, and the others non-$H$-edges, then since ${\rm Aut}(G)\subseteq{\rm Aut}(H)$, we can conclude that each automorphism of $G$ maps $H$-edges to $H$-edges and non-$H$-edges to non-$H$-edges. So assigning  each distinguishing edge labeling of $H$ to $G$ and assigning  non-$H$-edges a repeated label we make    a distinguishing edge labeling of $G$.  \qed

Since for two distinct simple non isomorphic, non rigid connected graphs, with no false twin  and no dominating vertex we have  ${\rm Aut}(G\star H) = {\rm Aut}(G\Box H)$, so a direct consequence of Lemmas \ref{nordspanning} and \ref{distindspann} is as follows:
\begin{theorem}\label{disindstrocartesi}
\begin{itemize}
\item[(i)] If $G$ and $H$ are two simple connected graphs, then $D'(G\star H)\leq D'(G\Box H)+1$.
\item[(ii)] If $G$ and $H$ are two  simple connected non isomorphic, non rigid  graphs with no false twin and no dominating vertex, then $D'(G\star H)\leq D'(G\Box H)$.
\end{itemize}
\end{theorem}

\begin{theorem}
	Let $ G$ be a  connected  graph with no false twin and no dominating vertex,	and $\star G^k$ the $k$-th power of $G$ with respect to the co-normal product. Then  for $k \geq 3$, $D'(\star G^k) = 2$. In particular, if $G$ is a rigid graph then  for $k \geq 2$, $D'(\star G^k) = 2$. 
\end{theorem} 
 \proof By Lemmas \ref{Lemma 2.3} and \ref{Lemma 4.1}, we can conclude that $G\star G$ has no false twin and no dominating vertex. We consider  the two following cases:
 
 Case 1)  Let $G$ be a non rigid graph. If $H= G\star G$, then $D(\star G^3) =2$ by  Theorem \ref{disindstrocartesi}(ii). Now by an induction on $k$, we have the result.
 
 Case 2) Let $G$ be a rigid graph. In this case, $|{\rm Aut}(G\star G)|=2$,  by Theorem \ref{thm2.17}, and so $D(G\star G) =2$. If $H:= G\star G$, then $|{\rm Aut}(G\star H)|=|{\rm Aut}(H)|$, by Theorem \ref{thm 2.15}. Hence $|{\rm Aut}(\star G^3)|=2$. By an induction on $k$ and using Theorem \ref{thm 2.15}, we obtain  $D(\star G^k) = 2$  for $k \geq 2$, where $G$ is a rigid graph.\qed
 
 \begin{theorem}\label{thm A}
 Let $G$ be a connected graph of order $n \geq 2$. Then $D'(G\star K_m)=2$ for every $m\geq 2$, except $D'(K_2\star K_2)=3$.
 \end{theorem}
 \proof Since  $|{\rm Aut}(G\star K_m)|\geq 2$, so $D'(G\geq K_m)=2$. With respect to the degree of vertices $G\star K_m$ we conclude that $G\star K_m$ is a traceable graph. We consider  the two following cases:
 
 Case 1) Suppose that  $n\geq 2$. If $m\geq 3$, or $m=2$, and $n\geq 4$, then the order of $G\star K_m$ is at least $7$, and so the result follows from Theorem \ref{traceble}. If $m=2$, $n=3$, then $G= P_3$ or $K_3$. In each case, it is easy to see that $D'(G\star K_m)=2$.
 
 Case 2) Suppose that  $n=2$. Then $G= K_2$, and so $G\star K_m= K_{2m}$. Thus $D'(G\star K_m)=2$ for $m\geq 3$, and  $D'(K_2\star K_2)=D'(K_4)=3$. \qed

By the value of the distinguishing index of Cartesian product of paths and cycles graphs in \cite{A. Gorzkowska R. Kalinowski and M. Pilsniak} and Theorem \ref{disindstrocartesi},  we can obtain this value for the co-normal product of them as the two following corollaries.

\begin{corollary}
\begin{itemize}
\item[\rm{(i)}] The co-normal product $P_m\star P_n$ of two paths of orders $m\geq2$ and $n\geq2$ has the distinguishing index equal to two, except $D'(P_2\star P_2)=3$.
\item[\rm{(ii)}] The co-normal product $C_m\star C_n$ of two cycles of orders $m\geq3$ and $n\geq3$ has the distinguishing index equal to two.
\item[\rm{(iii)}] The co-normal product $P_m\star C_n$  of orders $m\geq2$ and $n\geq3$ has the distinguishing index equal to two.
\end{itemize}
\end{corollary}
\proof
\begin{itemize}
	\item[\rm{(i)}]  If $n,m\geq 4$, then the result follows from  Theorem \ref{disindstrocartesi} (ii). If $n=2$ or $m=2$, then we have the result by Theorem \ref{thm A}. For the remaining cases,  with respect to the degree of vertices in $P_m\star P_n$ we obtain easily the distinguishing index.

	\item[\rm{(ii)}] If $n,m\geq 5$, then the result follows from  Theorem \ref{disindstrocartesi} (ii). If $n=3$ or $m=3$, then we have the result by Theorem \ref{thm A}. For the remaining cases we use of Hamiltonicity of $C_m\star C_n$ and Theorem \ref{traceble}.

	\item[\rm{(iii)}] If $n\geq 5$ and $m\geq 4$, then the result follows from  Theorem \ref{disindstrocartesi} (ii). If $n=3$ or $m=2$, then we have the result by Theorem \ref{thm A}. The remaining cases are $C_n\star P_3$ and $C_4\star P_m$. In the first case and with respect to the degree of vertices in $C_n\star P_3$ we obtain easily the distinguishing index. In the latter case, we use of Hamiltonicity of $C_4\star P_m$ and Theorem \ref{traceble}.
\qed
\end{itemize}

\end{document}